\documentclass[a4paper,reqno,11pt]{amsart}
\usepackage{amsmath,amssymb,amsfonts,amscd}

\title[Construction of $t$-structures and equivalences]{Construction of
$t$-structures and equivalences of derived categories}
\author[L. Alonso]{Leovigildo Alonso Tarr\'{\i}o}
\address[L. A. T.]{Departamento de \'Alxebra\\
Facultade de Matem\'a\-ticas\\
Universidade de Santiago de Compostela\\
E-15782  Santiago de Compostela, SPAIN}
\email{leoalonso@usc.es}

\author[A. Jerem\'{\i}as]{Ana Jerem\'{\i}as L\'opez}
\address[A. J. L.]{Departamento de \'Alxebra\\
Facultade de Matem\'a\-ticas\\
Universidade de Santiago de Compostela\\
E-15782  Santiago de Compostela, SPAIN}
\email{jeremias@usc.es}

\author[M. J. Souto]{Mar\'{\i}a~Jos\'e Souto Salorio}
\address[M. J. S. S.]{Facultade de Inform\'atica, Campus de Elvi\~na\\
Universidade da Coru\~{n}a\\
E-15071  A Coru\~{n}a, SPAIN}
\email{mariaj@udc.es}

\thanks{L.A.T. and A.J.L. partially supported by Spain's MCyT  and E.U.'s
FEDER research project BFM2001-3241, supplemented by Xunta de Galicia grant
PGDIT 01PX120701PR}

\subjclass[2000]{18E30 (primary); 14F05, 16D90 (secondary)}

\theoremstyle{plain}
\newtheorem{thm}{Theorem}[section]
\newtheorem{lem}[thm]{Lemma}
\newtheorem{cor}[thm]{Corollary}
\newtheorem{prop}[thm]{Proposition}
\newtheorem*{prob}{Problem}

\theoremstyle{remark}
\newtheorem*{rem}{Remark}

\theoremstyle{definition}

\newtheorem{cosa}[thm]{}

\newcommand{\A}{\mathcal{A}}

\newcommand{\CM}{{\mathcal C}}

\newcommand{\E}{{\mathcal E}}
\newcommand{\F}{{\mathcal F}}
\newcommand{\G}{{\mathcal G}}
\newcommand{\CH}{{\mathcal H}}
\newcommand{\CO}{{\mathcal O}}

\newcommand{\CL}{\mathcal{L}}
\newcommand{\CLA}{\CL_{\A}}

\newcommand{\V}{\mathcal{V}}
\newcommand{\U}{\mathcal{U}}
\newcommand{\SM}{{\mathcal S}}

\newcommand{\CC}{\boldsymbol{\mathsf{C}}}
\newcommand{\D}{\boldsymbol{\mathsf{D}}}
\newcommand{\K}{\boldsymbol{\mathsf{K}}}
\newcommand{\LL}{\boldsymbol{\mathsf{L}}}
\newcommand{\R}{\boldsymbol{\mathsf{R}}}
\newcommand{\T}{\boldsymbol{\mathsf{T}}}

\newcommand{\cc}{{\mathsf c}}
\newcommand{\bb}{{\mathsf b}}
\newcommand{\md}{\text{-}\mathsf{Mod}}
\newcommand{\pf}{\mathsf{cp}}
\newcommand{\qc}{\mathsf{qc}}

\newcommand{\PR}{\mathbf{P}}
\newcommand{\proj}{\PR^d_K}

\newcommand{\NN}{\mathbb{N}}
\newcommand{\ZZ}{\mathbb{Z}}

\newcommand{\holim}[1]{\begin{array}[t]{c} {\rm holim}\\[-7.5 pt]
 {\longrightarrow} \\[-7.5 pt] {\scriptstyle {#1}} \end{array}}
\newcommand{\dirlim}[1]{\begin{array}[t]{c} {\rm lim}\\[-7.5 pt]
 {\longrightarrow} \\[-7.5 pt] {\scriptstyle {#1}} \end{array}}

\newcommand{\lto}{\longrightarrow}

\newcommand{\noqed}{\renewcommand{\qed}{}}
\newcommand{\epi}{\twoheadrightarrow}

\newcommand{\inc}{\hookrightarrow}

\newcommand{\iso}{\tilde{\to}}
\newcommand{\liso}{\tilde{\lto}}
\newcommand{\imp}{\Rightarrow}

\DeclareMathOperator{\Hom}{Hom}
\DeclareMathOperator{\shom}{\CH\mathit{om}}
\DeclareMathOperator{\ext}{Ext}
\DeclareMathOperator{\sext}{\E\mathit{xt}}
\DeclareMathOperator{\rhom}{\R{}Hom}
\DeclareMathOperator{\Tot}{Tot}

\DeclareMathOperator{\Img}{Im}
\DeclareMathOperator{\cok}{Coker}
\DeclareMathOperator{\End}{End}
\DeclareMathOperator{\spec}{Spec}
\DeclareMathOperator{\real}{\mathsf{real}}
\DeclareMathOperator{\id}{id}
\DeclareMathOperator{\qco}{\mathsf{Qco}}

\newcommand{\ie}{{\it i.e.} }

\begin{document}

\begin{abstract} We associate a $t$-structure to a family of objects in
$\D(\A)$, the derived category of a Grothendieck category $\A$. Using 
general results on $t$-structures, we give a new proof of Rickard's
theorem on equivalence of bounded derived categories of modules. Also, we
extend this result to bounded derived categories of quasi-coherent sheaves on
separated divisorial schemes obtaining, in particular, Beilinson's
equivalences.
\end{abstract}

\maketitle

\section*{Introduction}

One important observation in modern homological algebra is that sometimes,
equivalences of derived categories do not come from equivalences of the
initial abelian categories. A remarkable example of this situation is that
of complexes of holonomic regular $\mathcal{D}$-modules and complexes of
locally constructible finite type sheaves of vector spaces over an analytic
variety. In fact, to a single regular holonomic $\mathcal{D}$-module
corresponds a full complex of constructible sheaves which are called
``perverse sheaves''. In turn, these sheaves are connected with intersection
cohomology that allows to recover Poincar\'e duality for singular analytic
spaces.

It is a crucial fact that a (bounded) derived category $\D^\bb$ may
contain as full subcategories abelian categories that are not isomorphic to
the initial one but whose (bounded) derived category is isomorphic to
$\D^\bb$. These facts have been systematized by Beilinson, Bernstein,
Deligne and Gabber in the work \cite{BBD} through the concept of
$t$-\emph{structure}. A $t$-structure in a triangulated category $\T$ is
given by functors which resemble the usual truncation functors in derived
categories and provide as a by-product an abelian full subcategory $\CM$ of
$\T$, its \emph{heart}, together with a homological functor on $\T$ with
values in $\CM$.

A derived category carries a natural $t$-structure but also can carry
others. It is clear that the construction of a $t$-structure without
recourse to an equivalence of derived categories is an important question.
In fact, this is how perverse sheaves are transposed to the \emph{\'etale}
context where it is not possible to define them through an equivalence.

But, up to now, there was no general method available to construct
$t$-structures, save for the glueing method of \cite[1.4]{BBD} whose main
application is to certain categories of sheaves over a site. In this paper,
we give such a general method for the derived category of a Grothendieck
abelian category (\ie an abelian category with a generator and exact direct
limits). These categories arise often in representation theory and
algebraic geometry. Our method parallels the fruitful method 
of localization used frequently in homotopy theory and algebra.

We believe our result is widely applicable. In this paper we will start
exploring these possibilities. We show that the remarkable theorem of
Rickard, that characterizes when two bounded derived categories of rings are
equivalent, can be deduced from our construction and standard facts on
$t$-structures. Moreover, we show that Rickard's result in \cite{Ric} can be
extended to certain schemes. Precisely, the bounded derived category of
quasi-coherent sheaves on a divisorial scheme is equivalent to the bounded
derived categories of modules over a certain ring if there exists a special
kind of ``tilting complex'' of sheaves. This generalizes Beilinson's
equivalence in \cite{B} and \cite{Bs}, in the sense that such equivalence
follows from this result.

Let us discuss now the contents of the paper. The first section recalls
the basic notions of $t$-structures and highlights the observation by
Keller and Vossieck that to give a $t$-structure $(\T^{\leq 0},
\T^{\geq 0})$ is enough to have the subcategory $\T^{\leq 0}$ and a right
adjoint for the canonical inclusion functor $\T^{\leq 0} \subset \T$. 
This has led us to seek for a proof in the spirit of Bousfield localization,
using techniques similar to those in \cite{AJS}. We will need to work in
the \emph{unbounded} derived category because the existence of coproducts is
essential for our arguments.

For a triangulated category $\T$, a suspended subcategory stable for
coproducts is called a cocomplete pre-aisle. In the second section we prove
that these subcategories are stable for homotopy direct limits. In the next
section we prove that if a cocomplete pre-aisle is generated by a set of
objects, then it is in fact a aisle, \ie a subcategory similar to 
$\T^{\leq 0}$. As in the case of Bousfield localizations, we proceed in two
steps. First, we deal with the case of derived categories of modules over a
ring, and next we treat the general case using the derived version of
Gabriel-Popescu embedding.

In the fourth section, we investigate the important question of whether the
$t$-structures previously defined restrict to the bounded derived category
of modules over a ring. Restriction to the upper-bounded category is easy,
however for the lower-bounded case some extra work is needed. We obtain this
restriction when the $t$-structure is generated by a compact object that also
generates the derived category in the sense of triangulated categories.

In the fifth section we deal with the analogous question in the context of
schemes. This section may be skipped if the reader is only interested in the
ring case. Trying to extend the arguments that work for rings in this
context, we need a condition on compatibility with tensor products in order
to localize to open subsets. This condition led us to define \emph{rigid
aisles} for which we give in Proposition \ref{rigeq} a criterion useful in
practice.

In the sixth section we see that the functor ``real'' defined in
\cite[Chap. III]{BBD} induces an equivalence of categories between the
bounded derived category of the heart of a $t$-structure and the bounded
derived category of modules over a ring precisely when the $t$-structure is
generated by a tilting object. The heart is equivalent to a category of rings
in this case ---recovering Rickard's equivalence theorem. The generalization
to the bounded derived category of quasi-coherent sheaves gives a similar
equivalence when the $t$-structure is rigid and generated by a tilting
object, thus making it equivalent to a bounded derived category of modules
over a ring. We explain how this result is related to Beilinson's
description of the bounded derived category of coherent sheaves on
projective space as a bounded derived category of modules over a finite
dimensional algebra.

We thank Amnon Neeman for discussions about these topics.

\section{Preliminaries on $t$-structures}

Let us begin fixing the notation and conventions.
Let $\T$ be a triangulated category whose translation autoequivalence is
denoted by $(-)[1]$ and its iterates by $(-)[n]$, with $n \in \ZZ$. A
$t$-\emph{structure} on $\T$ in the sense of Be\u{\i}linson, Bernstein,
Deligne and Gabber (\cite[D\'efinition~1.3.1]{BBD}) is a couple of full
subcategories $(\T^{\leq 0},\T^{\geq 0})$ such that, denoting
$\T^{\leq n} := \T^{\leq 0}[-n]$ and $\T^{\geq n} := \T^{\geq 0}[-n]$, the
following conditions hold:
  \begin{enumerate}
     \item[(t1)] For $X \in \T^{\leq 0}$ and $Y \in \T^{\geq 1}$,
                 $\Hom_{\T}(X,Y) = 0$.
     \item[(t2)] $\T^{\leq 0} \subset \T^{\leq 1}$ and 
                 $\T^{\geq 0} \supset \T^{\geq 1}$.
     \item[(t3)] For each $X \in \T$ there is a distinguished triangle
															  \[A \lto X \lto B \overset{+}{\lto}\] with 
															  $A \in \T^{\leq 0}$ and $B \in \T^{\geq 1}$.
  \end{enumerate}

The subcategory $\T^{\leq 0}$ is called the \emph{aisle} of the
$t$-structure, and $\T^{\geq 0}$ is called the \emph{co-aisle}. As usual for
a subcategory $\CM \subset \T$ we denote the associated orthogonal
subcategories as
$\CM^{\perp}= \{ Y \in \T / \Hom_{\T}(Z,Y) = 0, \;\forall Z \in \CM \}$
and 
${{}^{\perp}\CM} =  \{ Z \in \T / \Hom_{\T}(Z,Y) = 0, \;\forall Y \in \CM
\}$. Let us recall some immediate formal consequences of the definition.
\begin{prop} \label{first} 
Let $\T$ be a triangulated category, $(\T^{\leq 0},\T^{\geq 0})$ a
$t$-structure in $\T$, and $n \in \ZZ$, then
  \begin{enumerate}  
     \item $(\T^{\leq 0},\T^{\geq 1})$ is a pair of \emph{orthogonal}
           subcategories of $\T$, \ie $\T^{\geq 1} = {\T^{\leq 0}}^{\perp}$
           and $\T^{\leq 0} = {{}^{\perp}\T^{\geq 1}}$.
     \item The subcategories $\T^{\leq n}$ are stable for positive
           translations and the subcategories $\T^{\geq n}$ are
           stable for negative translations.
     \item The canonical inclusion $\T^{\leq n} \inc \T$ has a right adjoint
           denoted $\tau^{\leq n}$, and $\T^{\geq n} \inc \T$ a left adjoint
           denoted $\tau^{\geq n}$. Moreover, $X \in \T^{\leq n}$ if, and
           only if, $\tau^{\geq n+1}(X) = 0$, similarly for $\T^{\geq n}$. 
     \item For $X$ in $\T$ there is a distinguished triangle 
           \[\tau^{\leq 0}X \lto X \lto \tau^{\geq 1}X \overset{+}{\lto}.\]
     \item The subcategories $\T^{\leq n}$ and  $\T^{\geq n}$ are stable
           for extensions, \ie given a distinguished
           triangle $X \to Y \to Z \overset{+}{\to}$, if $X$ and $Z$ belong
           to one of these categories, so does $Y$.     
  \end{enumerate}
\end{prop}

\begin{proof}
The statement $(i)$ follows easily from the axioms and $(ii)$ a restatement
of (t2). The statements $(iii)$ and $(iv)$ are proved in \cite[Prop.
1.3.3]{BBD}. Finally, let's sketch a proof of $(v)$ for $\T^{\leq n}$, the
one for $\T^{\geq n}$ is dual. By translation, we can even restrict to 
$\T^{\leq 0}$. Apply $\Hom_{\T}(-,\tau^{\geq 1}Y)$ to the triangle and get
the long exact sequence of abelian groups:
\[
\dots \to \Hom_{\T}(Z,\tau^{\geq 1}Y) \to \Hom_{\T}(Y,\tau^{\geq 1}Y)
\to \Hom_{\T}(X,\tau^{\geq 1}Y) \to \dots
\]
where $X$ and $Z$ belong to $\T^{\leq 0} = {{}^{\perp}\T^{\geq 1}}$.
Therefore both extreme homs are zero and by exactness
\[\Hom_{\T}(\tau^{\geq 1}Y,\tau^{\geq 1}Y) \cong \Hom_{\T}(Y,\tau^{\geq
1}Y) = 0\] which means that $\tau^{\leq 0}Y = Y$, \ie $Y \in \T^{\leq 0}$.
\end{proof}

The interest of this notion lies in the fact that a $t$-structure 
$(\T^{\leq 0},\T^{\geq 0})$ in a triangulated category $\T$ provides a full
abelian subcategory of $\T$, its \emph{heart}, defined as 
$\CM :=  \T^{\leq 0} \cap \T^{\geq 0}$ and a homological functor, namely
$H^0 := \tau^{\geq 0} \tau^{\leq 0}$,  with values in $\CM$, \emph{cfr.}
\cite[Thm. 1.3.6]{BBD}. One would like to know natural conditions that
guarantee the existence of a $t$-structure in a triangulated category.  
We will show how to transpose homotopical techniques to get means of
constructing $t$-structures.

A fundamental observation is to use the structure of the subcategories
$\T^{\leq n}$ and  $\T^{\geq n}$. It is clear that, in general, they are not
triangulated subcategories but they come close. In fact, each subcategory
$\T^{\leq n}$ has the structure of \emph{suspended category} in the sense of
Keller and Vossieck \cite{KVsus}. Let us recall this definition.

An additive category $\U$ is suspended if and only if is graded by an
additive translation functor $T$ (sometimes called \emph{shifting}) and
there is class of diagrams of the form $X \to Y \to Z \to TX$ (often denoted
simply 
$X \to Y \to Z \overset{+}{\to}$) called \emph{distinguished triangles} such
that the following axioms, analogous to those for triangulated categories in
Verdier's \cite[p.~266]{vtc} hold:
\begin{enumerate}  
     \item[(SP1)] Every triangle isomorphic to a distinguished one is
                  distinguished. For $X \in \U$, 
                  $0 \to X \overset{\id}{\to} X \to 0$ is a distinguished
                  triangle. Every morphism $u:X \to Y$ can be completed to a
                  distinguished triangle 
                  $X \overset{u}{\to} Y \to Z \to TX$
     \item[(SP2)] If $X \overset{u}{\to} Y \to Z \to TX$ is a distinguished
                  triangle in $\U$ then so is 
                  $Y \to Z \to TX \overset{Tu}{\to} TY$.
     \item[(SP3)] = (TR3) in Verdier's \emph{loc.~cit}.
     \item[(SP4)] = (TR4) in Verdier's \emph{loc.~cit}.
\end{enumerate}
 
The main difference with triangulated categories is that the translation
functor in a suspended category may not have an inverse and therefore some
objects can not be shifted back. The formulation of axioms (SP1) and (SP2)
reflect this fact.

In this paper, we will only consider suspended subcategories of a
triangulated category $\T$ and we implicitly assume that they are full and
that the translation functor is the same of $\T$. If 
$(\T^{\leq 0},\T^{\geq 0})$ is a $t$-structure, the aisle
$\T^{\leq 0}$ is a suspended subcategory of $\T$ whose distinguished
triangles are diagrams in $\T^{\leq 0}$ that are distinguished triangles
in $\T$ (Proposition \ref{first}). Moreover, the aisle
$\T^{\leq 0}$ determines the $t$-structure because the co-aisle
$\T^{\geq 0}$ is recovered as $(\T^{\leq 0})^{\perp}[1]$. The terminology
``aisle'' and ``co-aisle'' comes from \cite{KVais}.

The following observation that sent us towards the right direction for our
objective is the following result, due to Keller and Vossieck:
\begin{thm} (\cite[Section 1]{KVais}) A suspended subcategory $\U$ of a
triangulated category $\T$ is an aisle (\ie $(\U,\U^{\perp}[1])$ is a
$t$-structure on $\T$) if and only if the canonical inclusion functor $\U
\inc \T$ has a right adjoint.
\end{thm}


Our arguments, however, are logically independent of this fact. The
knowledgeable reader will recognize a localization situation: to construct a
$t$-structure, one basically looks for a right adjoint to the inclusion of a
subcategory. In the case of triangulated localizations a basic concept is
\emph{localizing} subcategories: those triangulated subcategories stable
for coproducts. We will use the analogous concept in the context of
$t$-structures. 

To respect Keller and Vossieck's terminology we will call
\emph{pre-aisle} a suspended subcategory $\U$ of a triangulated category $\T$
where the triangulation in $\U$ is given by the triangles which are
distinguished in $\T$ and the shift functor is induced by the one in $\T$.
We see easily that to check that a full subcategory $\U$ of $\T$ is a
pre-aisle, it is enough to verify that
\begin{itemize}
   \item For any $X$ in $\U$, $X[1]$ is also in $\U$.
   \item Given a distinguished triangle $X \to Y \to Z \to X[1]$, if $X$
         and $Z$ belong to $\U$, then so does $Y$.
\end{itemize}
Once these two facts hold for $\U$, the verification of axioms
(SP1)--(SP4) is immediate. 

Let $\U$ be an aisle, \ie a suspended subcategory of $\T$ such that
$(\U,\U^{\perp}[1])$ is a $t$-structure. We will denote by $\tau_\U^{\leq n}$
and $\tau_\U^{\geq n}$ ($n \in \ZZ$) the truncation functors associated to
this $t$-structure, \ie the corresponding adjoints to the canonical
inclusions $\U[-n] \inc \T$ and $\U^{\perp}[1-n] \inc \T$, respectively.

\begin{lem}
Let $\T$ be a triangulated category and let $\U$ be an aisle.
If a family of objects of  $\U$ has a coproduct in $\T$, then
it already belongs to $\U$.
\end{lem}
\begin{proof}  For $Y$ in $\T$,
\[\Hom(\oplus X_i, \tau_\U^{\geq 1} Y) = 
\prod \Hom(X_i, \tau_\U^{\geq 1}
Y);\] so if $\tau_\U^{\geq 1} X_i =0$ for every $i$, then 
$\tau_\U^{\geq 1}(\oplus X_i)= 0$, so $\oplus X_i$ is in $\U$.
\end{proof}

The following fact was observed first in the case of localizing subcategories
by B\"okstedt and Neeman in \cite{BN}.

\begin{cor} \label{dirsum}
A direct summand of an object in a aisle $\U$ belongs to $\U$.
\end{cor}
\begin{proof}
Use the argument in \emph{loc. cit.}\/ or, alternatively, Eilenberg's
swindle, as in \cite[footnote, p.~227]{AJS}.
\noqed\end{proof}

If $\U$ is a pre-aisle in $\T$ stable for (arbitrary) coproducts in $\T$,
then it will be called \emph{cocomplete}. From the previous discussion, we
see that in triangulated categories in which coproducts exist, aisles are
cocomplete pre-aisles. Specifically, we try to address the following:

\begin{prob}\label{prob}
For a triangulated category $\T$, give conditions that
guarantee that a cocomplete pre-aisle is an aisle.
\end{prob}

In this paper, we will solve this problem in two cases. Let $\A$ be
a Grothendieck category, that means $\A$ is an abelian category with a
generator $U$ and exact (set-indexed) filtered direct  limits. We will give a
sufficient condition for a cocomplete pre-aisle to be an aisle for the
derived category of $\A$. This will be that there is a set (as opposed to a
class) of objects that ``generates'' the pre-aisle.   To fix notation,
denote by $\CC(\A)$ the category of complexes of objects of $\A$ (also a
Grothendieck category). Denote by $\K(\A)$ and $\D(\A)$ the homotopy and
derived categories of $\A$, with their usual structure of triangulated
categories. The proof will rely in a construction by homotopy limits of
diagrams made by maps in $\CC(\A)$ rather that in $\D(\A)$, so it can not be
transposed to a general triangulated category. However if the set that
``generates'' the cocomplete pre-aisle is made of \emph{compact} objects
(see below) then by a similar but simpler argument (possibly well-known),
the analogous result is obtained for any triangulated category. We will
treat this case in an appendix.

\section{Cocomplete pre-aisles and homotopy colimits}

In this section we will show that homotopy colimits of objects in a
cocomplete pre-aisle $\U$ of $\D(\A)$ belong to $\U$.  We will follow the
strategy in \cite{AJS} for localizing subcategories, so the proofs will be
only sketched as briefly as possible.

For convenience of the reader, let us recall the construction of the
homotopy colimit of a diagram of complexes as used in \cite{AJS}. Let $G =
\{G_s, \mu_{s\,t}\:/\: s \leq t \in \Gamma\}$ be a filtered directed system
of $\CC(\A)$. Let $r > 0$ and $s \in \Gamma$, $W_s^r$ is the set of chains in
the ordered set $\Gamma$ of length $r$ which start in $s$. Define a
bicomplex $B(G)$ by: $B(G)^{k\,j} := 0$ if $k > 0$, $B(G)^{0\,j} :=
\bigoplus_{s \in \Gamma}G^j_s$ and
$B(G)^{k\,j} := \bigoplus_{s \in \Gamma,w \in W_s^{-k}}G^j_{s,w}$ if 
$k < 0$, where by $G^j_{s,w}$  we denote $G^j_s$ indexed by a (fixed)
chain $w$ of $W_s^{-k}$. If $k < 0$, we will denote a map from the generator 
$x:U \to G^j_{s,w} \to B(G)^{k\,j}$ by $(x;s < s_1 < \cdots < s_{-k})$
where $w$ is the chain $s < s_1 < \cdots < s_{-k}$ and then the horizontal
differential $d_1$ is defined by the following formula:
\begin{small}
\begin{equation*}
   \begin{split}
&d_1^{k\,j}(x;s < s_1 < \cdots < s_{-k})   := \\
& = (\mu_{s\,s_1}(x); s_1 < \cdots < s_{-k}) + \sum_{i=1}^{-k}((-1)^i x;s <
\cdots < \check{s_i} < \cdots < s_{-k}),
   \end{split}
\end{equation*}
\end{small}%
where the symbol $\check{s_i}$ means that $s_i$ is suppressed from
the chain. If $k \geq 0$ then $d_1$ is 0. 
The differential $d_2$ is induced by the one of the complexes $G_{s}$.
The {\it homotopy direct limit} of the system is defined as the totalization
(by coproducts) of the bicomplex $B(G)$, and we denote it by
$\holim{s\in\Gamma}G_s$.

\begin{lem} \label{milnor}
Let $G = \{G_n / n \in \NN\}$ be a directed system in $\CC(\A)$. Let
$\U$ be a cocomplete pre-aisle of  $\D(\A)$.  If $G_n \in \U$, for every
$n\in \NN$, then \[\dirlim{n \in \NN}G_n \in \U,\]
where the limit is taken in $\CC(\A)$.
\end{lem}

\begin{proof}
Consider the Milnor distinguished triangle in $\D(\A)$ (\cite[3.3]{AJS})
\[
\bigoplus_{n \in \NN}G_n \overset{1-\mu}{\lto} \bigoplus_{n \in \NN}G_n 
\lto \dirlim{n \in \NN}G_n \overset{+}{\lto}.
\]%
Being $\U$ cocomplete, every $G_n$ is in $\U$, so the first two 
objects in this triangle are in $\U$ and therefore the third, because a
pre-aisle is stable for mapping cones.
\end{proof}

\begin{rem}
For a directed system $\{G_n \:/\: n \in \NN \}$ in a triangulated category
$\T$, there is a construction of a homotopy direct limit by B\"okstedt and
Neeman that makes sense for any triangulated category (see \cite{BN}). It is
defined by the triangle
\[
\bigoplus_{n \in \NN}G_n \overset{1-\mu}{\lto} \bigoplus_{n \in \NN}G_n 
\lto C(1-\mu) \overset{+}{\lto},
\]%
where $C(1-\mu)$ denotes the mapping cone of $1-\mu$. It is clear that given
a cocomplete pre-aisle $\U$ in $\T$, if every $G_n \in \U$, then \( C(1-\mu)
\in \U.\)
Moreover, by \cite[Theorem 2.2]{AJS}, we have the isomorphism in $\D(\A)$
\[\holim{n \in \NN}G_n \:\iso \dirlim{n \in \NN}G_n.\]
Also is is clear that 
\[C(1-\mu) \:\iso \dirlim{n \in \NN}G_n.\]
Thus, it will do no harm to identify these three objects if one is working
inside $\D(\A)$.
\end{rem}

\begin{prop} \label{limte}
Let $\U$ be a cocomplete pre-aisle of $\D(\A)$.
Let $\Gamma$ be a filtered ordered set and 
$G = \{G_s, \mu_{s\,t}\:/\: s \leq t \in \Gamma\}$ be a
directed system in $\CC(\A)$ such that $G_s$ belongs to $\U$ for every $s
\in \Gamma$. Then \[\holim{s\in \Gamma}G_s \] also belongs to $\U$.
\end{prop}

\begin{proof}
By definition, \[\holim{s\in \Gamma}G_s  = \Tot(B(G)^{\cdot\,\cdot}),\]
where
$B(G)^{\cdot\,\cdot}$ is a bicomplex whose columns $B(G)^{i\,\cdot}$ are
coproducts of the objects in the system, so they belong to $\U$. By
\cite[Lemma 3.2]{AJS}, there is sequence $F_n$ of subcomplexes of
$\Tot(B(G)^{\cdot\,\cdot})$ such that 
\[\Tot(B(G)^{\cdot\,\cdot}) = \dirlim{n \in \NN}F_n\]
in $\CC(\A)$.
So by the previous lemma we only have to show that every  $F_n$
is in $\U$ for every $n\in\NN$. Indeed, $F_0 = B(G)^{0\,\cdot}$, a coproduct
of objects of $\U$. For $n > 0$, take $\nu_n$ a graded splitting of the
inclusion $F_{n-1} \subset B(G)^{-n+1\,\cdot}[n-1]$. The composition
\[ B^{-n\,\cdot}[n-1] \xrightarrow{d_1[n-1]} B^{-n+1\,\cdot}[n-1]
\overset{\nu_n}{\lto} F_{n-1}^{\cdot} 
\]%
defines a map of complexes $h_n$ and it can be completed to a triangle
\[
B^{-n\,\cdot}[n-1]  \overset{h_n}{\lto}
F_{n-1}^{\cdot} {\lto} F_{n}^{\cdot} \overset{+}{\lto}.
\]
So, $F_n$ is the cone of a map from a \emph{positive} shift of a column of
$B(G)^{\cdot\,\cdot}$ to $F_{n-1}$, which belongs to $\U$ by induction,
 therefore $F_n \in \U$ by Lemma \ref{milnor}.
\end{proof}

\section{Construction of $t$-structures}
Let $R$ be a ring. In this paper, we will only consider associative rings
with unit. Denote by $R\md$ the category of left $R$-modules. As usual, we
abbreviate $\CC(R\md)$, $\K(R\md)$ and $\D(R\md)$, by $\CC(R)$, $\K(R)$ and
$\D(R)$, respectively. We denote by  $Q_R:\K(R) \to \D(R)$ the canonical
functor.

\begin{lem} \label{lete}
Let $\SM = \{E_{\alpha} / \alpha \in A \}$ be a set of objects of a
triangulated category $\T$ with coproducts. Let $\U$ be the smallest
cocomplete pre-aisle that contains the objects in $\SM$. The groups
$\Hom_{\T}(E_{\alpha}[j],B) = 0$ for every $j \geq 0$ and $\alpha \in A$
if, and only if, $B$ is in $\U^{\perp}$.
\end{lem}

\begin{proof}
The ``if'' part is trivial. Let us prove the ``only if'' part. Suppose that
$B \in \T$ is such that  $\Hom_{\T}(E_{\alpha}[j],B) = 0$ for every $j \geq
0$ and $\alpha \in A$, we have to show that  $\Hom_{\T}(X,B) = 0$ for every
$X \in \U$. Let $\V$ be the full subcategory of $\U$ whose objects $X$ are
such that $\Hom_{\T}(X[j],B) = 0$, for every $j \geq 0$. $\V$ is clearly
stable for positive shiftings and for extensions. Finally, $\V$ is stable for
coproducts, because
\[
\Hom_{\T}(\bigoplus_{i \in I}X_i[j],B) \simeq
\prod_{i \in I}\Hom_{\T}(X_i[j],B)
\]%
Now, $\SM$ is contained in $\V$, therefore $\V = \U$.
\end{proof}

\begin{prop} \label{tconsDR}
Let $\U$ the smallest cocomplete pre-aisle of $\D(R)$ which contains an
object $E \in \CC(R)$. Then $\U$ is an aisle in $\D(R)$.
\end{prop}

\begin{proof} 
First, we can assume that the complex $E$ is K-projective by the well-known 
existence of resolutions (see, for instance, \cite[Proposition 4.3]{AJS}).
Therefore, for $k \in \ZZ$ every $E[k]$ is also  K-projective and maps from
$E[k]$ in the derived category are represented by actual maps of complexes. 
Observe that the conditions (t1) and (t2) of the definition of $t$-structure
hold automatically for $(\U,\U^{\perp}[1])$, therefore it is enough to prove
(t3), \ie to construct for any $M \in \CC(R)$ a distinguished triangle $N
\to M \to B
\overset{+}{\to}$ with $N \in \U$ and $B \in \U^{\perp}$. 

The construction of $B$ is by a transfinite induction and is parallel to the
construction given in the proof of \cite[Proposition, 4.5]{AJS}.
Let $\gamma$ be the least ordinal such that
\( \#(\gamma) > \#\cup_{p \in \ZZ}E^p \) 
and let $I$ be a set of ordinals that contains $\gamma$.

Let 0 be the minimum of $I$, we define $B_0 := M$.

If $s \in I$ has a predecessor $s - 1$, suppose by induction that $B_{s-1}$
is already constructed. Take $\Omega_{s-1} := \bigcup _{k \in \NN}
\Hom_{\CC(R)}(E[k],B_{s-1})$, and let  
$\alpha_{s-1} : \bigoplus_{\Omega_{s-1}}E[k] \to B_{s-1}$ be  given by the
universal property of the coproduct. We define 
$\mu_{s-1\,s} : B_{s-1} \to B_s$ by the canonical distinguished triangle:
\[ \begin{CD}
\bigoplus_{\Omega_{s-1}}E[k] @>\alpha_{s-1}>> B_{s-1} @>\mu_{s-1\,s}>> B_s 
@>>> \left( \bigoplus_{\Omega_{s-1}}E[k] \right) [1]
\end{CD} ,\]%
where $B_s$ is the mapping cone of $\alpha_{s-1}$.
For any $i < s$, let $\mu_{i\,s} := \mu_{s-1\,s} \circ \mu_{i\,s-1}$. 

If $s \in I$ has no predecessor, take 
\[
B_s := \dirlim{i < s}B_i 
\text{ and for } i < s \text{ , } \mu_{i\,s} = \dirlim{i < j < s}\mu_{i\,j},
\]%
both limits constructed as objects and maps in $\CC(R)$.

We have got a directed a system in $\CC(R)$,
$\{ B_s, \mu_{s\,t} / s \leq t \in I\}$, such that every $\mu_{s\,t}$ is
semi-split and if $f \in \Hom_{\CC(R)}(E[k],B_s)$, where $s < t$ and $k \in
\NN$, then $\mu_{s\,t} \circ f$ is  homotopic to zero.

We have to show that $B := B_{\gamma}$ belongs to $\U^{\perp}$. By the
previous lemma ($\SM = \{ E \}$) it is enough to show that
$\Hom_{\D(R)}(E[k],B) =  0$, for each $k \in \NN$. Let $g : E[k] \to B$ be a
map of complexes, we will show that $g$ is homotopic to zero. 
Let $B'_{i} := B_i \cap \Img(g)$.
The ordinal $\gamma$ is a cardinal, so it has no predecessor and
\[
\dirlim{i < \gamma}B'_i = \Img(g)
\]%
We claim that there is an index $s_0 < \gamma$ such that 
\[\dirlim{i < s_0}B'_i = B'_{s_0}.\]
Thus, there is a $s_0 \in I$ such that  $g$ factors through
$\mu_{s_0\,\gamma}$, and by the properties of our directed system $g$ is
homotopic to zero.

Let us show the claim now. Arguing by contradiction, suppose that for
every $s \in I$, there is a $t > s$ such that the map $B'_s \inc B'_t$ is
not an epimorphism. Define sets  $J_s := \{r \geq s / B'_r \inc B'_{r+1}
\; \text{is not an epimorphism}\}$. Each $J_s$ is not empty, so $J := J_0$
is a cofinal subset of $I$. But $\gamma$ is a regular ordinal by
Hausdorff's theorem (\cite[3.11 Proposition, p. 135]{Le}) and therefore,
$\#(J) =\#(\gamma)$.

On the other hand, we can define a map
 \[ \phi : J \to \cup_{p \in \ZZ}\Img(g)^p\] by 
$\phi(s) = \alpha_s$ where $\alpha_s \in \cup_{p \in \ZZ}{B'}_{s+1}^p$
but $\alpha_s \notin \cup_{p \in \ZZ}{B'}_s^p$. The map $\phi$ is
injective and $\Img(g)$ is a quotient of $E[k]$, therefore,
\(\#(J) \leq \#(\cup_{p \in \ZZ}\Img(g)^p) < \#(\gamma), \) a
contradiction.

The only thing left to prove is that in the triangle 
$N \to M \to B \overset{+}{\to}$ we have that $N \in \U$. For each 
$i < \gamma$, it is clear that defining $N_i$ by the distinguished triangle
\begin{equation} \label{trii}
N_i \lto M \overset{\mu_{0\,i}}{\lto} B_i \lto N_i[1]
\end{equation}
we can take $N_i[1] = \cok(\mu_{0\,i})$ (in $\CC(R)$) and that
we have $N = \dirlim{i < \gamma}N_i$, where for $s \leq t$ the transition
maps $N_s \to N_t$ of this system are induced by $\mu_{s\,t}:B_s \to B_t $.
We will check that $N_i \in \U$ for every $i \in I$ and this will finish
the proof.

First, $N_0 = 0$, so it belongs to $\U$. If $i$ has a predecessor,
consider the triangles like (\ref{trii}) whose second maps are
$\mu_{0\,{i-1}}$ and $\mu_{0\,i}$. Consider also the distinguished triangle:
\[
\oplus_{\Omega_{i-1}}E[k] \lto B_{i-1} \overset{\mu_{i-1\,i}}{\lto} B_{i}
\lto \left( \oplus_{\Omega_{i-1}}E[k]\right)[1].
\]%
We have $\mu_{i-1\,i} \circ \mu_{0\,{i-1}} = \mu_{0\,i}$. we can apply the
octahedral axiom which gives us a distinguished triangle
\[
N_{i-1} {\lto} N_{i}  \lto
 \oplus_{\Omega_{i-1}}E[k] \lto N_{i-1}[1],
\]%
where $N_{i-1}$ and $\oplus_{\Omega_{i-1}}E[k]$ belong to $\U$.
This implies $N_i$ is also in $\U$. Finally, if $i$ has no predecessor,
$N_i$ is a direct limit of complexes $N_s$ that belong to $\U$ by induction. 

Applying \cite[Theorem 2.2]{AJL} we have that
\[\holim{s < i}N_s \simeq \dirlim{s < i}N_s\]
and by Proposition \ref{limte}, the homotopy limit belongs to $\U$, so we
reach the desired conclusion.
\end{proof}

Given an aisle $\U$ defined as the smallest that contains an object $E$, we
will denote the truncation functors associated to the $t$-structure
$(\U,\U^{\perp}[1])$ by $\tau_E^{\leq n} =\tau_\U^{\leq n}$ and
$\tau_E^{\geq n} = \tau_\U^{\geq n}$ ($n \in \ZZ$).

\begin{rem} The hypothesis of the previous result is extended immediately
to several generators, the smallest cocomplete pre-isle that contains a
\emph{set} of objects agrees with the one that contains its coproduct
because being cocomplete implies being stable for direct summands (Corollary
\ref{dirsum}).
\end{rem}

\begin{cosa}\label{gab-pop} Our next step is to generalize the  previous
result to the derived category of a Grothendieck category $\A$.
By Gabriel-Popescu embedding (\cite{GP}, see also \cite[Ch. X, Theorem
4.1]{St}),  $\A$ is a quotient (in the sense of Abelian
categories) of $R\md$ by a thick subcategory whose objects are the
\emph{torsion objects} of a hereditary torsion theory, where $R$ is
$\End_{\A}(U)$ and $U$ denotes a generator of $\A$. This means there is a
couple of functors
\begin{equation} \label{adjunction}
R\md \underset{i}{\overset{a}{\rightleftarrows}} \A
\end{equation}
where $a$ is an exact functor and $i$ is left-exact full, faithful and
right adjoint to $a$. The torsion objects are those $R$-modules that
are sent to zero by $a$.

By the existence of Bousfield localization in the derived category of a
ring (\cite[Proposition 4.5]{AJS}), the situation extends to derived
categories. Precisely, there is a diagram of triangulated functors:
\begin{equation}
\D(R) \underset{\mathbf{i}}{\overset{\mathbf{a}}{\rightleftarrows}} \D(\A).
\end{equation}
The functor $a$ is exact, therefore it has a canonical extension to the
derived category, which we denote $\mathbf{a}$. This functor has a right
adjoint $\mathbf{i}$ which extends $i$ whose existence is not immediate but
follows from \cite[Proposition 5.1]{AJS}. This says that there is a
localizing subcategory $\CLA$ of $\D(R)$ generated by a set of complexes
such that the quotient $\D(R)/\CLA$ is equivalent to $\D(\A)$ and the
associated localization functor is identified with $\mathbf{i} \circ
\mathbf{a}$. Note that $\mathbf{a} \circ \mathbf{i} = \id_{\D(\A)}$. Let 
$F$ be a generator of the localizing subcategory $\CLA$. Using Lemma
\ref{lete} we can check that the smallest cocomplete pre-aisle that contains
the set of objects $\{F[-n] / n \in \NN\}$ is in fact $\CLA$. This is an
essential point for the proof of the next result.
\end{cosa}

\begin{thm} \label{tcons}
Let $\U$ the smallest cocomplete pre-aisle of $\D(\A)$ which contains an
object $E \in \CC(\A)$. Then $\U$ is an aisle in $\D(\A)$.
\end{thm}
\begin{proof}
Keeping the notation of the previous discussion let $\U'$ the smallest
cocomplete pre-aisle in $\D(R)$ that contains $\mathbf{i}E$ and $F[-n]$ for
all $n \in \NN$.
Then by
Proposition \ref{tconsDR} it is an aisle. Consequently for an object $M \in
\D(\A)$, there is a distinguished triangle in $\D(R)$,
\[N_0 \lto \mathbf{i}(M) \lto B_0 \overset{+}{\lto},
\] where $N_0 \in \U'$ and $B_0 \in {\U'}^{\perp} \subset \CLA^{\perp}$.
This gives a triangle in $\D(\A)$:
\[ \mathbf{a}(N_0) \lto M \lto  \mathbf{a}(B_0) \overset{+}{\lto}.
\]
The pre-aisle $\U$ is the essential image by $\mathbf{a}$ of $\U'$, so
it is clear that $\mathbf{a}(N_0)$ belongs to  $\U$. Let us check that
$\mathbf{a}(B_0)$ is in $\U^{\perp}$. Let $j \geq 0$, using the fact that
$E$ comes from $\D(\A)$ and $B_0$ is $\CLA$-local, we see that
\[ \Hom_{\D(\A)}(E[j], \mathbf{a}B_0) 
\simeq \Hom_{\D(R)}(\mathbf{i}E[j],\mathbf{i}\mathbf{a}B_0)) = 0
\]
Applying Lemma \ref{lete} again, we conclude.
\end{proof}

\section{Boundedness and $t$-structures}

From here on, we will apply our method of construction of
$t$-structures to the problem of characterizing equivalences of bounded
derived categories. In this section and the next we will deal with the issue
of when a $t$-structure defined on the full category $\D(\A)$, restricts to
a subcategory defined through a boundedness condition.

\begin{prop}\label{demenos} Let $E \in \D^-(\A)$. The $t$-structure
defined by $E$ restricts to $\D^-(\A)$, in other words, the associated
truncation functors, $\tau_E^{\geq n}$ and $\tau_E^{\leq n}$ ($n \in
\ZZ$) take upper bounded objects to upper bounded objects.
\end{prop}

\begin{proof} It is enough to check the assertion for $\tau_E^{\leq 0}$ and
$\tau_E^{\geq 1}$. Note that the aisle $\U$ generated by $E$ is such that
$\U \subset \D^-(\A)$, because the operations that allow us to construct
$\U$ from $E$ ---translation and extension--- preserve this condition,
therefore, for every $M \in \D^-(\A)$, $\tau_E^{\leq 0}(M) \in \U \subset
\D^-(\A)$. Now considering the triangle,
\[\tau_E^{\leq 0}M \lto M \lto \tau_E^{\geq 1}M \overset{+}{\lto}
\]
$M$ and $\tau_E^{\leq 0}M$ are in $\D^-(\A)$, therefore, 
$\tau_E^{\geq 1}M \in \D^-(\A)$.
\end{proof}

Our next task will be to get conditions for the restriction of a
$t$-structure to $\D^+(\A)$. In the rest of this section we will treat the
case of derived categories of modules.

\begin{cosa}\label{compact} We need a standard definition in the theory of
triangulated categories. Let $\T$ be a triangulated category. An object $E$
of $\T$ is called \emph{compact} if the functor $\Hom_{\T}(E,-)$ commutes
with arbitrary coproducts. Another way of expressing the condition is that a
map from $E$ to a coproduct factors through a finite subcoproduct.
\end{cosa}

The following is Rickard's criterion for compact objects in derived
categories over rings. See \cite[Proposition 6.3 and its proof]{Ric}. Using
results of Neeman, we are able to give a simpler proof.

Let $\T$ be a triangulated category. A triangulated subcategory $\SM
\subset \T$ is called \emph{thick} if a direct summand of an object of $\SM$,
also belongs to $\SM$. Observe that a localizing subcategory is thick,
because being closed for coproducts implies being closed for direct summands,
as follows form Corollary \ref{dirsum}.

\begin{lem} \label{critR}
Let $R$ be a ring. In $\D(R)$ the compact objects are those isomorphic to
bounded complex of finite-type projective modules.
\end{lem}

\begin{proof} The triangulated category $\D(R)$ is generated by $R$ that is
trivially compact and bounded. We apply Neeman's result (\cite[Lemma
2.2]{Ntty}) that says that the smallest thick subcategory of a compactly
generated triangulated category that contains all the translations of a
compact generator is the full subcategory whose objects are the compact
ones. Therefore the full subcategory of compact objects is the smallest
thick subcategory of $\D(R)$ that contains $R$. Now the full subcategory of
bounded complexes of projectives is thick because it is stable for triangles
and direct summands, and all of its objects are compact as can be
checked easily, therefore it agrees with the subcategory of compact objects.
\end{proof}

\begin{cosa}In a triangulated category $\T$, a family of objects 
$\{E_{\alpha}\}_{\alpha \in A}$ is called a \emph{generating family} if for
every $M \in \T$, $M \neq 0$ there is a $j \in \ZZ$ and an $\alpha \in A$
such that $\Hom_{\T}(E_{\alpha}[j],M) \neq 0$. If the family is formed by a
single object, then this object is usually called a generator.
\end{cosa}

The following result will be key for applications.

\begin{prop}\label{demas} Let $R$ be a ring. Let $E$ be a compact generator
of $\D(R)$, and denote by $\U$ the aisle generated by $E$. Then every object 
$M \in \U^{\perp}$ belongs to  $\D^+(R)$.
\end{prop}

\begin{proof}
Fix $M \in \U^{\perp}$. Let $\SM_M$ be the full subcategory of $\D(R)$ whose
objects are those $N$ in $\D(R)$ such that $\rhom^{\cdot}_R(N,M) \in
\D^+(\ZZ)$. 

The subcategory $\SM_M$ is a thick subcategory of $\D(R)$.
Indeed, it is clearly triangulated and a direct summand $N'$ of an object
$N$ of $\SM_M$, also belongs to $\SM_M$ because
$\rhom^{\cdot}_R(N',M)$ is a direct summand of $\rhom^{\cdot}_R(N,M)$.
The object $E$ and all of its translations $E[k]$ ($k \in \ZZ$) belong to
$\SM_M$ because $M \in \U^{\perp}$ which implies that for all $j \leq 0$
\[ H^j(\rhom^{\cdot}_R(E,M)) = \Hom_{\D(R)}(E[-j],M) = 0,\]
so $\rhom^{\cdot}_R(E,M) \in \D^{\geq 1}(\ZZ)$.

We have shown that $\SM_M$ is a thick subcategory of $\D(R)$ that contains
all translations of a compact generator. We apply again Neeman's result
(\cite[Lemma 2.2.]{Ntty}). It follows that all the compact objects are
contained in $\SM_M$, in particular $R \in \SM_M$.
But this means that $\rhom^{\cdot}_R(R,M) \in \D^+(\ZZ)$, in other
words, $H^j(M) = H^j(\rhom^{\cdot}_R(R,M)) = 0$, for $j \ll 0$, \ie $M \in
\D^+(R)$, as desired.
\end{proof}

\begin{thm} \label{boundring}
Let $R$ be a ring and $E$ a compact generator of $\D(R)$. The
$t$-structure defined in $\D(R)$ by\/ $\U$, the smallest aisle that contains
$E$,  restricts to a $t$-structure on $\D^\bb(R)$.
\end{thm}

\begin{proof} 
Consider again the triangle associated to $M \in \D(R)$
\[\tau_E^{\leq 0}M \lto M \lto \tau_E^{\geq 1}M \overset{+}{\lto}.
\]
By Lemma \ref{critR} we can assume that $E$ is bounded, in particular, $E \in
\D^-(R)$ and, by Proposition \ref{demenos}, the $t$-structure defined by $E$
restricts to $\D^-(R)$. This means that when $M \in \D^\bb(R)$, then both 
$\tau_E^{\leq 0}M$ and $\tau_E^{\geq 1}M \in \D^-(R)$.

On the other hand, by Proposition \ref{demas}, we always have that 
$\tau_E^{\geq 1}M \in \D^+(R)$, therefore if $M \in \D^\bb(R)$,
$\tau_E^{\leq 0}M \in \D^+(R)$ also. Putting all this together,
$\tau_E^{\leq 0}M$ and $\tau_E^{\geq 1}M \in \D^\bb(R)$.
\end{proof}

\section{Boundedness for quasi coherent sheaves}

The previous result can be extended for sheaves on separated divisorial
schemes. Recall that a scheme $(X,\CO_X)$ is called divisorial if it is
quasi-compact and quasi-separated and has an ample family of line bundles,
see \cite{I}. An illuminating example is given by a quasi-projective variety
(over a field, say). 
\medskip

Denote by $\qco(X)$ the category of quasi-coherent sheaves on a scheme $X$.
If $X$ is separated, the category $\D(\qco(X))$ has an internal tensor
functor $-\otimes^{\LL}_{\CO_X}-$ that can be computed by quasi-coherent
flat resolutions in either variable, see \cite[1.1]{AJL}. 

Let $d \in \NN$, an aisle $\U$ in $\D(\qco(X))$ is called $d$-\emph{rigid} if
for every  $\G \in \U$ and $\F \in \D^{\leq 0}(\qco(X))$ we have that
$\F\otimes^{\LL}_{\CO_X}\G \in \U[-d]$. If $\U$ is $d$-\emph{rigid} for some
$d$, we will say simply that $\U$ is rigid. The next result gives a useful
characterization of rigid aisles in $\D(\qco(X))$ when $X$ is a divisorial
scheme.

\begin{prop} \label{rigeq} Let $(X,\CO_X)$ be a divisorial separated scheme.
Let $\U$ be an aisle of $\D(\qco(X))$ and $d \in \NN$. The following are
equivalent:
  \begin{enumerate}
     \item $\U$ is $d$-rigid.
     \item For every $\G \in \U$ and $\F \in \qco(X)$ we have
           that $\F\otimes^{\LL}_{\CO_X}\G \in \U[-d]$.
     \item For every $\G \in \U$ and $\F$ a locally free sheaf in 
           $\qco(X)$ we have that $\F\otimes_{\CO_X}\G =
           \F\otimes^{\LL}_{\CO_X}\G \in \U[-d]$.
     \item For every $\G \in \U$, every invertible sheaf $\CL$ from the ample
           family of sheaves in $\qco(X)$, there is a $k_0 \in \NN$ such
           that $\CL^{\otimes k}\otimes_{\CO_X}\G \in \U[-d]$ if $k < k_0$.
  \end{enumerate}
\end{prop}

\begin{proof} The assertions $(i) \imp (ii) \imp (iii) \imp (iv)$ are
trivial. To prove $(iv) \imp (i)$ take a complex $\F \in \D^{\leq
0}(\qco(X))$. The complex  $\F\otimes^{\LL}_{\CO_X}\G$ can be computed by
means of a flat resolution $\E \to \F$, where each component $\E^i$ is an
arbitrary coproduct of sheaves $\CL^{\otimes k}$ where $\CL$ is
an invertible sheaf from the ample family and $k \ll 0$. Therefore we have
the quasi-isomorphism  $\E\otimes_{\CO_X}\G \to \F\otimes^{\LL}_{\CO_X}\G$.
The complex $\E\otimes_{\CO_X}\G$ is defined as the totalization of the
bicomplex $\E^{\cdot}\otimes_{\CO_X}\G^{\cdot}$. Fix $i \leq 0$, the complex
$\E^{i}\otimes_{\CO_X}\G^{\cdot}$ belongs to $\U[-d]$, by $(iv)$ and the
previous observation using that an aisle is stable for coproducts and the
derived tensor product of sheaves commutes with coproducts. This allows us
to apply the same argument as in Proposition \ref{limte}. We define a
filtration $\CH_n$ of $\E\otimes_{\CO_X}\G$ such that 
$\CH_0 = \E^{0}\otimes_{\CO_X}\G^{\cdot}$. For $n > 0$, there is a map
$\nu_n$ from which we define $\CH_n$ inductively by the triangle:
\[(\E^{-n}\otimes_{\CO_X}\G^{\cdot})[n-1] \overset{\nu_n}{\lto} 
\CH_{n-1}^{\cdot}  {\lto} 
\CH_{n}^{\cdot}    \overset{+}{\lto}. \]
Thus every $\CH_n \in \U[-d]$ by induction. Clearly, in $\CC(\qco(X))$ 
\[\dirlim{n}\CH_n = \E\otimes_{\CO_X}\G,\] 
and it follows that $\E\otimes_{\CO_X}\G \in \U[-d]$ by Lemma \ref{milnor}.
\end{proof}

\begin{cor} In the previous situation, if $X$ is affine, then every aisle in
$\D(\qco(X))$ is rigid.
\end{cor}

\begin{proof}In this case $\{\CO_X\}$ is an ample family.
\end{proof}


\begin{cosa} Let $(X,\CO_X)$ be a scheme. Let $\CL$ be an invertible sheaf
on $X$ and $s : \CO_X \to \CL$ a global section. For any $n \in \NN$, there
is a map $\CL^{\otimes n} \to \CL^{\otimes n+1}$ defined by $t \mapsto
t \otimes s$. These maps form a directed system of $\CO_X$-modules whose
limit we will denote
\[\CO_X[s^{-1}] := \dirlim{n}\CL^{\otimes n}.
\]
The notation is explained by the fact that if $X$ is affine, say 
$X = \spec(R)$, and $s \in R$ is identified with a global section of the
sheaf $\CO_X = \widetilde{R}$, then $\CO_X[s^{-1}] = \widetilde{R_s}$.
\end{cosa}

\begin{lem} \label{Riso}
Let $(X,\CO_X)$ be a divisorial separated scheme and $\CL$ be an invertible
sheaf from the given ample family. Let $X_s$ be the open set whose
complementary set is the set of zeros of a global section  $s : \CO_X \to
\CL$. Denote by $j : X_s \inc X$ the canonical inclusion map. Let $\F \in
\D(\qco(X))$, then we have an isomorphism (in $\D(\qco(X))$):
\[ \R{}j_*j^*\F \liso \CO_X[s^{-1}]\otimes_{\CO_X}\F.
\]
\end{lem}

\begin{proof} By \cite[Lemme 2.2.3.1]{I} the morphism 
$j : X_s \inc X$ is affine, therefore $j_*: \qco(X_s) \to \qco(X)$ is exact.
Being an open embedding, $j$ is flat, so $j^*$ is exact, too. Therefore, we
have the isomorphism $\R{}j_*j^*\F \iso j_*j^*\F$ for every $\F \in
\D(\qco(X))$. So, it is enough to treat the case in which  $\F$ is a single
quasi-coherent sheaf. In fact, this is a restatement of \cite[6.8.1]{ega1}.

Consider the canonical commutative diagram
\[ \begin{CD}
\F        @>>>       \CO_X[s^{-1}] \otimes_{\CO_X} \F \\
@VVV                   @V{\beta}VV                           \\          
j_*j^*\F  @>{\alpha}>> \CO_X[s^{-1}] \otimes_{\CO_X} j_*j^*\F
\end{CD} \]
The maps $\alpha$ and $\beta$ are isomorphisms. Indeed, consider a covering
of $X$ made up by affine trivialization open subsets, \ie every $U =
\spec(R)$ in the covering is such that $\CL|_U \iso \widetilde{R}$, and the
section $s|_U$ is identified with a certain  $f \in R$. Suppose also that 
$\F|_U \iso \widetilde{M}$ for a certain $R$-module $M$. Observe that the
open set $X_s \cap U$ is identified with $\spec{R_f}$. With this
identifications $\alpha|_U$ corresponds to the canonical isomorphism 
$\widetilde{M_f} \iso \widetilde{R_f} \otimes_{\CO_X|_U} \widetilde{M_f}$ 
and $\beta|_U$ corresponds to  $\widetilde{R_f} \otimes_{\CO_X|_U}
\widetilde{M} \iso \widetilde{R_f} \otimes_{\CO_X|_U} \widetilde{M_f}$ 
\end{proof}

For divisorial schemes there is a characterization of compact objects in
$\D_{\qc}(X)$, the category of complexes of sheaves of $\CO_X$-modules with
quasi-coherent homology, due to Thomason and Trobaugh, see \cite[Theorem
2.4.3.]{tt}. On the other hand, it is well-known that for a quasi-compact
separated scheme we have an equivalence $\D_{\qc}(X) \simeq \D(\qco(X))$,
see, for instance \cite[Corollary 5.5]{BN} or, for a different proof,
\cite[Proposition 1.3]{AJL}. We give here a proof of this criterion for
$\D(\qco(X))$.

\begin{lem} \label{critTT}
Let $(X,\CO_X)$ be a divisorial, separated scheme. In $\D(\qco(X))$ the
compact objects are those isomorphic to bounded complexes of locally free,
finite type sheaves.
\end{lem}

\begin{proof} We apply again Neeman's result (\cite[Lemma 2.2.]{Ntty}) to
$\D(\qco(X))$ and an ample family of invertible sheaves as
compact generators. Then we mimic the procedure of Lemma \ref{critR}.
\end{proof}

\begin{prop}\label{demasqc} Let $(X,\CO_X)$ be a divisorial, separated
scheme. Let $\E$ be a compact generator of $\D(\qco(X))$, and denote by $\U$
the aisle generated by $\E$. Suppose that $\U$ is rigid. Then every object 
$\F \in \U^{\perp}$ belongs to  $\D^+(\qco(X))$.
\end{prop}

\begin{proof} As in the case of Proposition \ref{demas} we only need to check
that \[ \CH^q(\F) = \sext_{\CO_X}^q(\CO_X, \F)\]
vanishes for $q \ll 0$. But this is a local question, therefore
we can cover $X$ by a finite family of open affines,  
$U_{\lambda} = \spec(R_{\lambda})$, where  $\lambda \in \{1, \dots, n \}$.
Denote by $j_{\lambda} : U_{\lambda} \inc X$ the canonical inclusion. We
choose these open affines as complementary sets of zeros of sections
$s_{\lambda} : \CO_X \to \CL_{\lambda}$ of invertible sheaves
$\CL_{\lambda}$ from the ample family: $U_{\lambda}$  will be the
locus where $s_{\lambda}$ does not vanish. The complex $\E$ restricted to
each $\spec(R_{\lambda})$ can be represented by a bounded complex of free
$R_{\lambda}$-modules (after refining the covering, if necessary). Define
complexes $E_{\lambda}$ and $F_{\lambda}$ of $R_{\lambda}$ modules by
$\widetilde{E_{\lambda}} = \E|_{\spec(R_{\lambda})}$ and 
$\widetilde{F_{\lambda}} = \F|_{\spec(R_{\lambda})}$. 

The complex $E_{\lambda}$ is a compact generator in $\D(R_{\lambda})$.
Indeed, by Lemma \ref{critTT} the covering by open affines can be chosen in
such a way that $E_{\lambda}$ are bounded complexes of projective modules of
finite type which are compact as objects of $\D(R_{\lambda})$. Let us check
that $E_{\lambda}$ is a generator of $\D(R_{\lambda})$. For a non zero
object $M$ of $\D(R_{\lambda})$, its derived extension 
$\R{}j_{\lambda *}\widetilde{M}$ is again non zero because
$j^*_{\lambda}\R{}j_{\lambda *}\widetilde{M} = \widetilde{M} \neq 0$.
Therefore, there is a non zero map 
$\E[r] \to \R{}j_{\lambda *}\widetilde{M}$, for some $r \in \ZZ$ which, by
adjunction, gives a non zero map $E_{\lambda}[r] \to M$.

Suppose that $\U$ is $d$-rigid. Let $\U_{\lambda}$ be the aisle of
$\D(R_{\lambda})$ generated by $E_{\lambda}[d]$. We will see that
$F_{\lambda} \in \U_{\lambda}^{\perp}$. We have to check that  
\[\Hom_{\D(R_{\lambda})}(E_{\lambda}[r], F_{\lambda}) = 0\] for every 
$r \geq d$. First, using once more \cite[Theorem 2.2]{AJS} we have that
\[
\holim{n}\CL^{\otimes n} \iso \dirlim{n}\CL^{\otimes n} = \CO_X[s^{-1}]
\]
Now we have the following chain of isomorphisms
\begin{equation*}
   \begin{split}
        \Hom_{\D(R_{\lambda})}(E_{\lambda}[r], F_{\lambda}) 
              & \iso \Hom_{\D(\qco(U_{\lambda}))}(j_{\lambda}^*\E[r],
                      j_{\lambda}^*\F) \\
              & \iso \Hom_{\D(\qco(X))}(\E[r],
                     \R{}j_{\lambda *}j_{\lambda}^*\F) \\
              & \iso \Hom_{\D(\qco(X))}(\E[r], \CO_X[s_{\lambda}^{-1}] 
                      \otimes \F)\\
              & \iso \Hom_{\D(\qco(X))}(\E[r], 
                     (\!\!\holim{n}\CL^{\otimes n}) \otimes \F) \\
              & \iso \dirlim{n} \Hom_{\D(\qco(X))}(\E[r],
                     \CL^{\otimes n} \otimes\F) \\
              & \iso \dirlim{n} \Hom_{\D(\qco(X))}(\E[r] \otimes
                     \CL^{\otimes -n}, \F) \\
              & = 0,\\
   \end{split}
\end{equation*}
where the first isomorphism holds because $U_{\lambda}$ is affine,  the third
using Lemma \ref{Riso}, the fourth, following the previous remark, the
fifth, using that $\E$ is a compact object in  $\D(\qco(X))$ and the rest
are the usual adjunction maps. The last equality follows from the fact that
$\U$ is $d$-rigid.

Now, given $\lambda \in \{1, \dots, n \}$, Proposition \ref{demas} gives us a
bound on $q$ for the vanishing of $\ext_R^q(R_{\lambda}, F_{\lambda})$.
As these bounds are in finite quantity, we take the minimum and conclude.
\end{proof}

\begin{thm} \label{boundsch}
Let X be a separated divisorial scheme $X$. If $\E$ is compact and 
generator in $\D(\qco(X))$, the $t$-structure defined by $\E$ in
$\D(\qco(X))$ restricts to a $t$-structure on $\D^\bb(\qco(X))$.
\end{thm}

\begin{proof}  Mimic the proof of Theorem \ref{boundring}, using Proposition
\ref{demasqc} instead of Proposition \ref{demas} and that 
$\E \in \D^-(\qco(X))$ by Lemma \ref{critTT}.
\end{proof}

\section{Derived equivalences and $t$-structures}

First, we need the following easy Lemma.

\begin{lem} \label{critperp}
Let $\T$ be a triangulated category and $\U$ be the smallest cocomplete
pre-aisle in $\T$ that contains an object $E$. If $E$ is \emph{compact} then
$\Hom_{\T}(E,X[j]) = 0$ for all  $j > 0$ and all  $X \in \U$ if, and only
if, $\Hom_{\T}(E,E[j]) = 0$ for all
$j > 0$.
\end{lem} 

\begin{proof}
We will use again an idea from a previous Lemma. Consider the full
subcategory $\SM$ whose objects are those $X \in \U$ such that
$\Hom_{\T}(E,X[j]) = 0$ for every $j > 0$. Clearly, $\SM \subset \U$ is
suspended and stable for coproducts ($E$ is compact), therefore if 
$E \in \SM$, necessarily $\SM = \U$. 
\end{proof}

\begin{cosa}
In a triangulated category $\T$, an object $X$ is called \emph{exceptional}
if, and only if, \[\Hom_{\T}(X,X[j]) = 0 \text{, for all } j \neq 0.\]
E. g., an object $M$ in $\D(R)$ ($R$ a ring) is exceptional if  
$\ext_{R}^j(M,M) = 0$ for all $j \neq 0$. In the literature, an exceptional
compact generator is often called a \emph{tilting} object.
\end{cosa}

\begin{prop}\label{progen} Let $\A$ be a Grothendieck category.
Let $\U$ be the smallest aisle in $\D(\A)$ which contains an object
$E$. If $E$ is compact and exceptional, then it belongs to the heart 
$\CM = \U^{\perp}[1] \cap \U$ of the $t$-structure $(\U,\U^{\perp}[1])$. The
category $\CM$ is abelian, cocomplete and $E$, as an object of $\CM$, is
\emph{small}, \emph{projective} and \emph{generator}.
\end{prop}

\begin{proof}
First, $E$ belongs to $\U$ and to $\U^{\perp}[1]$ because it is
exceptional, therefore $E \in \CM$. 
The category $\CM$ is abelian by \cite[Th\'eor\`eme 1.3.6]{BBD}. To see that
it is cocomplete it is enough to see that it has coproducts. This is true
for $\U$ and to $\U^{\perp}[1]$ by Proposition \ref{aico} and therefore also
for $\CM$.

Take now any object $Y$ in $\CM$. It belongs to $\U$ by definition. If we
have also that $\Hom_{\CM}(E,Y) = \Hom_{\D(\A)}(E,Y) = 0$ then it already
belongs to $\U^{\perp}$, therefore to $\U \cap \U^{\perp}$, so $Y = 0$, this
shows that $E$ is a generator in $\CM$ (in the sense of abelian categories).

Let  $0 \to X \to Y \to Z \to 0$ be an exact sequence in $\CM$. It
comes from a triangle in $\D(\A)$, namely,
\(X \lto Y \lto Z \overset{+}{\lto}.\)
As $Z \in \U^{\perp}[1]$, we have that $Z[-1] \in \U^{\perp}$. But then
$\Hom_{\D(\A)}(E,Z[-1])=0$. The complex $E$ is exceptional, then by Lemma
\ref{critperp}, $\Hom_{\D(\A)}(E,X[1])=0$ because $X$ belongs to $\U$. Taking
these facts into account, we apply $\Hom_{\CM}(E,-) = \Hom_{\D(\A)}(E,-)$ to
the triangle and obtain an exact sequence:
\[0 \to \Hom_{\CM}(E,X) \lto \Hom_{\CM}(E,Y) \lto \Hom_{\CM}(E,Z) \to 0,\]
therefore $E$ is a projective object of $\CM$. The object $E$ is compact in
$\D(\A)$, thus it is small in $\CM$ because the functor $\Hom_{\CM}(E,-)$
commutes with coproducts.
\end{proof}

For the next results, we need some notation. Let $\T$ be a triangulated
category and $\U$ an aisle in $\T$. The bounded part of $\T$ with respect to
$\U$ is the full subcategory defined by:
\[
\T^\bb_\U := \bigcup_{a,b \in \ZZ} \U[a] \cap \U^{\perp}[b]
\]
The following Lemma is similar to \cite[Proposition
3.1.16]{BBD}.

\begin{lem} \label{bbdual}
Let $(\U, \U^{\perp}[1])$ be a $t$-structure in $\D(\A)$ and $\CM$ its
heart. The functor $\real : \D^\bb(\CM) \to \D(\A)^\bb_\U$, defined
in \cite[3.1.10]{BBD}, is an equivalence of categories if, and only if, 
for every objects $A$ and $B$ in $\CM$, $ n > 0$ and $f \in
\Hom_{\D(\A)}(A,B[n])$, there is an epimorphism $A' \epi A$ in $\CM$ that
erases $f$. 
\end{lem}

\begin{proof} To check that $\real$ is fully faithful, we modify slightly the
first part of the proof in \emph{loc.~cit.} The functor $\real$ induces a map
\[\Hom_{\D^\bb(\CM)}(A,B[n]) \to \Hom_{\D(\A)}(A,B[n])\]
that is an isomorphism for $n = 0$.
For $n > 0$, the functors $\Hom_{\D^\bb(\CM)}(-,B[n])$ are Yoneda's Ext in
$\CM$, therefore make part of a (contravariant) universal $\delta$-functor.
And the same is true for $\Hom_{\D(\A)}(-,B[n])$, because the conditions
of Grothendieck's characterization of universal $\delta$-functors
(\cite[Proposition 2.2.1]{toh}) hold by hypothesis.

The essential image of $\real$ is $\D(\A)^\bb_\U$ by the same argument as in
\cite[Proposition 3.1.16]{BBD}.
\end{proof}

For convenience, for a separated divisorial scheme $X$, let us say that a
tilting complex of quasi-coherent sheaves $\E$ is \emph{special} if the
aisle generated by $\E$ in $\D(\qco(X))$ is rigid.

\begin{prop} \label{core}
Let the category $\A$ be either $R\md$ for an arbitrary ring $R$ or $\qco(X)$
for a separated divisorial scheme $X$. Let $E$ be a tilting object of
$\D(\A)$, which we suppose special in the scheme case. Let $\CM$ be the
heart of the $t$-structure that $E$ defines in $\D(\A)$. Then,  $\D^\bb(\CM)$
is equivalent to $\D^\bb(\A)$
\end{prop}

\begin{proof} Let $\U$ denote the aisle determined by $E$. The object $E \in
\D^-(\A)$, because it is bounded (Lemma \ref{critR} and Lemma \ref{critTT}),
therefore, we have that $\U[a] \subset \D^-(\A)$ for every $a \in \ZZ$.
Using Proposition \ref{demas} for $\A = R\md$ and Proposition
\ref{demasqc} for  $\A = \qco(X)$ we see that $\U^{\perp}[b] \subset
\D^+(\A)$ for every $b \in \ZZ$. As a consequence $\D(\A)^\bb_\U \subset
\D^\bb(\A)$.

Let us show that $\D(\A)^\bb_\U = \D^\bb(\A)$. Given $M \in \D^\bb(\A)$, we
will find $a, b \in \ZZ$ such that  $M \in \U[a]$ and $M \in \U^{\perp}[b]$.
Let us check first that $\rhom^{\cdot}(E,M)$ is a bounded complex. Indeed,
in the case that  $\A = R\md$ the compact object $E$ can be represented by a
bounded complex of projective modules (Lemma \ref{critR}) and therefore
$\rhom^{\cdot}(E,M) \iso \Hom^{\cdot}(E,M)$ which is clearly bounded. In
case $\A = \qco(X)$, the compact complex $E$ can be represented by a bounded
complex of locally free finite-type modules (Lemma \ref{critTT}) and we have
that
\(\R\shom^{\cdot}(E,M) \iso \shom^{\cdot}(E,M)\), because locally free
modules are acyclic for $\shom^{\cdot}(-,M)$ from which it follows that
\(\R\shom^{\cdot}(E,M)\) is bounded. Now, $X$ being separated, the derived functor $\R\Gamma(X,-)$ can
be computed via \v{C}ech resolutions (\cite[1.4.1]{ega3}) and these
resolutions are bounded because $X$ is quasi-compact. Therefore
$\R\Gamma(X,-)$ takes bounded complexes in $\A$ to bounded complexes of
abelian groups and
\[\rhom^{\cdot}(E,M) \cong \R\Gamma(X,\R\shom^{\cdot}(E,M))
\] is a bounded complex as claimed. Then there exist $a,b \in \ZZ$ such that
$\Hom_{\D(\A)}(E[j],M) = H^{-j}\rhom^{\cdot}(E,M) = 0$ for all $j \geq b$ and
$j < a$.

Note that $\U^{\perp}[b] = \U[b]^{\perp}$ and the aisle $\U[b]$ is generated
by $E[b]$. We apply Lemma \ref{lete} and conclude that $M \in \U^{\perp}[b]$
because $\Hom_{\D(\A)}(E[j],M) = 0$ for all $j \geq b$. 

To see that $M \in \U[a]$, we prove that in the distinguished triangle
\[\tau^{\leq -a}_{\U}M \lto M \lto \tau^{\geq 1-a}_{\U}M
\overset{+}{\lto},\] the object
$B := \tau^{\geq 1-a}_{\U}M$ is $0$. The aisle $\U[a]$ is generated by
$E[a]$, that is a tilting complex, and then
\[\Hom_{\D(\A)}(E[a+j],E[a]) = \Hom_{\D(\A)}(E[j],E) = 0,\] for
$j < 0$. By Lemma \ref{critperp}, if $N \in \U[a]$ we have
that $\Hom_{\D(\A)}(E[a+j],N) = 0$ for all $j<0$. Also,
$\Hom_{\D(\A)}(E[a+j],M) = 0$ for all $j<0$. Applying the cohomological
functor $\Hom_{\D(\A)}(E,-)$ to the distinguished triangle above we obtain
an exact sequence from which it follows that $\Hom_{\D(\A)}(E[a+j],B) = 0$
for $j < 0$.  On the other hand, we have that 
$\Hom_{\D(\A)}(E[a+j],B) = 0$ for $j \geq 0$ because $B \in \U[a]^\perp$.
Summing up, as $E$ is a generator, $B = 0$.

Finally, consider the functor $\real : \D^\bb(\CM) \to \D(\A)^\bb_\U$. To
see that it is an equivalence, it is only left to check that the hypothesis
of the previous Lemma hold. Let $A$ and $B$ in $\CM$, $ n > 0$ and $f \in
\Hom(A,B[n])$. By Proposition \ref{progen}, $E$ is projective and a generator
in $\CM$. We can take an epimorphism $A' \epi A$, where $A'$ is a direct sum
of copies of $E$. And this epimorphism erases $f$, as required.
\end{proof}

\begin{thm} \label{eqR}Let $R$ and $S$ be rings. We have an equivalence:
   \begin{enumerate}
       \item The categories $\D^\bb(S)$ and $\D^\bb(R)$ are equivalent.
       \item There is a tilting object $E$ in $\D(R)$, such that
             $\Hom_{\D(R)}(E,E) = S$
   \end{enumerate}
\end{thm}

\begin{proof} $(i) \Rightarrow (ii)$ The object of $\D^\bb(R)$ corresponding
to $S$ by the equivalence is clearly tilting.

$(ii) \Rightarrow (i)$ Consider the $t$-structure generated by the tilting
object $E$ in $\D(R)$. Let $\CM$ be its heart. By Proposition \ref{core},
$\D^\bb(R)$ is equivalent to $\D^\bb(\CM)$. Also, by Proposition
\ref{progen}, the object $E$ belongs to $\CM$ and it is a small projective
generator. But an abelian category with a small projective generator is a
module category by \cite[Chap. IV, Theorem 4.1 p. 104]{Mit}. Then the  ring
$S := \Hom_{\CM}(E,E)$ is such that $S\md$ is equivalent to
$\CM$ and consequently, $\D^\bb(S)$ is equivalent to $\D^\bb(R)$.
\end{proof}

\begin{rem} The previous theorem gives another proof of (b) $\Leftrightarrow$
(e) in \cite[Theorem~6.4.]{Ric}. We have to say that Rickard was also
inspired in the construction given in \cite{BBD} to get a totalization of
bicomplexes ``up to homotopy''. The difference is our explicit use of
the functor $\real$. This shows that constructing a $t$-structure
is a general means for establishing an equivalence.

We also point out that some technicalities from sections 4 and 5
could be avoided if we had an analogous of the functor $\real$ for
unbounded derived categories. This construction would be feasible once we
had a satisfactory theory of totalization for unbounded bicomplexes which,
at present, does not seem obvious to get.
\end{rem}

\begin{thm}\label{tsch} Let $X$ be a separated divisorial scheme. If there
is a special tilting object $\E$ in $\D(\qco(X))$, then the ring 
$R= \Hom_{\D(\qco(X))}(\E,\E)$ is such that $\D^\bb(\qco(X))$ is equivalent
to $\D^\bb(R)$.
\end{thm}

\begin{proof} It follows from Proposition \ref{core}, using again the
characterization of module categories in \cite{Mit}
as in the previous Theorem.
\end{proof}

\begin{cosa} Let us see now how Beilinson's equivalence can be interpreted
under the light of the previous results. Let $\proj$ be the $d$-dimensional
projective space over a field $K$. Denote by $\CO(1)$ the canonical ample
invertible sheaf and by $\Omega^1_\PR$ the sheaf of differential forms.
Consider the families of locally free sheaves $\SM_\Lambda := \{\CO_\PR,
\Omega^1_\PR(1), \dots, \Omega^d_\PR(d)\}$ and 
$\SM_S : = \{\CO_\PR, \CO_\PR(-1), \dots, \CO_\PR(-d)\}$. 
They are obviously made of compact objects in
$\D(\qco(\proj))$. From \cite[Lemma 2]{B}, the objects
\begin{align*} 
\E_\Lambda & := 
   \CO_\PR \oplus \Omega^1_\PR(1) \oplus \dots \oplus  \Omega^d_\PR(d) \\ 
\E_S & := 
   \CO_\PR \oplus \CO_\PR(-1) \oplus \dots \oplus  \CO_\PR(-d) 
\end{align*}
are exceptional. Also, from the proof of \cite[Theorem,
p. 215]{B} they are generators of $\D(\qco(\proj))$. Let
\begin{align*}   
R_\Lambda &:= \End_{\D(\qco(\proj))}(\E_\Lambda) \\ 
R_S       &:= \End_{\D(\qco(\proj))}(\E_S).
\end{align*} 
Denote simply by $R$  either $R_\Lambda$ or $R_S$, and, analogously, let
$\E$ be either $\E_\Lambda$ or $\E_S$, respectively. For a triangulated
category $\T$, denote by $\T_\pf$ the full subcategory of compact objects in
$\T$. We have the following.
\end{cosa}

\begin{cor} There is an equivalence between the categories
$\D^\bb(\qco(\proj))$ and $\D^\bb(R)$. Moreover, this equivalence induces one between
$\D^\bb_\cc(\qco(\proj))$ and $\D(R)_\pf$.
\end{cor} 

\begin{proof} In view of Theorem \ref{tsch}, the only thing left to check is
that $\E$ is special in either case, in other words, that the aisle $\U$
generated by $\E$ (or equivalently, by the families $\SM_\Lambda$ or
$\SM_S$) is rigid.

The key argument in Beilinson's papers is that one can express every coherent
sheaf $\F$ over $\proj$ as the 0-th homology of a complex concentrated in
degrees from $-d$ to $d$ (denoted $L^\cdot_\Lambda$ and $L^\cdot_S$ in
\cite[Lemma 3]{Bs}) whose components are formed by coproducts of the objects
in the families $\SM_\Lambda$ and $\SM_S$, respectively. These resolutions
are called the \emph{Beilinson monads} of the sheaf $\F$. It follows that
$\F$ can be obtained by a process of producing mapping cones, coproducts and
\emph{positive} translations up to $d$ of objects in the families and
therefore $\F \in \U[-d]$. 

It follows that, for $k \ll 0$ and $i \in \{0, \dots d\}$, the
objects $\Omega^i_\PR(i+k)$, twists of the generators of the aisle
$\U_\Lambda$, belong to $\U_\Lambda[-d]$; and analogously for the sheaves
$\CO(-i+k)$ and the aisle $\U_S$. Now it is easy to see that, once this fact
holds for the generators of an aisle $\U$, it holds also for all the objects
in $\U$. Therefore we can apply Proposition \ref{rigeq} to $\U_\Lambda$
or $\U_S$, respectively, and conclude that they are both rigid.

Finally, the object $\E$
is projective and generator viewed as $R$-module, and therefore is a compact
generator of $\D(R)$. By \cite[Lemma 2.2]{Ntty} $\D(R)_\pf$ is the smallest
thick subcategory of $\D^\bb(R)$ that contains $\E$. The functor $\real$
extends the inclusion of the heart to a map between derived categories, so
it takes $\E$ to itself. But $\E$ is also a compact generator
of $\D(\qco(\proj))$ and by the cited Lemma, the smallest thick
subcategory of $\D^\bb(\qco(\proj))$ containing $\E$
is $\D(\qco(\proj))_\pf$. The functor $\real$ preserves thick subcategories
and therefore both are equivalent. Finally, $\D^\bb_\cc(\qco(\proj)) =
\D(\qco(\proj))_\pf$ by Auslander--Buchsbaum--Serre theorem, because $\proj$
is a regular scheme.
\end{proof}


\renewcommand{\thesection}{\Alph{section}}
\setcounter{section}{1}
\setcounter{thm}{0}

\section*{Appendix: Compactly generated $t$-structures}

In this section we will give a solution to the \textbf{Problem} in
page \pageref{prob} of Section 1 for general triangulated categories under
the hypothesis that the generating objects for the complete pre-aisle are
compact. This gives a criterion for the existence of $t$-structures
applicable in very general settings. Also, the good behavior of these kind
of objects will allow us to get information about the $t$-structure already
in the case of a derived category, as treated in the main text.

\begin{thm}\label{tcompact}
Let $\SM = \{E_{\alpha} / \alpha \in A \}$ a \emph{set} of compact objects
in a triangulated category $\T$. Let $\U$ the smallest cocomplete pre-aisle
of $\T$ which contains the family $\SM$. Then $\U$ is an aisle in $\T$.
\end{thm}

\begin{proof}
This argument is inspired in \cite[Lemma 1.7, p. 554]{Ntty}.
Let $M \in \T$, let us show how to construct a distinguished triangle
$N \to M \to B \overset{+}{\to}$ with $N \in \U$ and $B \in \U^{\perp}$. 

This time ordinary induction will do. Let $B_0 := M$. Suppose $B_{n-1}$ is
already constructed. Let $\Omega_{n-1} := \bigcup _{\alpha \in A, k \in
\NN}  \Hom_{\CC(R)}(E_{\alpha}[k],B_{n-1})$ and define
$B_{n}$ by the triangle
\begin{equation}\label{trinum}
\bigoplus_{\Omega_{n-1}}E_{\alpha}[k] \overset{\rho}{\lto} B_{n-1}
\overset{\mu}{\lto} B_n \overset{+}{\lto},
\end{equation}%
where $\rho$ is defined by the universal property of the coproduct, as in
Proposition \ref{tconsDR}. Let $B = \holim{n \in \NN} B_n$ which, by the
remark after Lemma \ref{milnor}, may be defined by the distinguished triangle
\begin{equation*}
\bigoplus_{n \in \NN}B_n \overset{1-\mu}{\lto} \bigoplus_{n \in
\NN}B_n  \lto B \overset{+}{\lto}.
\end{equation*}%
The object $B$ belongs to $\U^{\perp}$. Indeed, using Lemma \ref{lete}
it is enough to check that for any $\alpha \in A$ and $k \geq 0$,
$\Hom_{\T}(E_{\alpha}[k],B) = 0$. But, applying the functor
$\Hom_{\T}(E_{\alpha}[k],-)$ to the previous triangle,
\[
\Hom_{\T}(E_{\alpha}[k],\holim{n \in \NN} B_n) \cong
\dirlim{n \in \NN} \Hom_{\T}(E_{\alpha}[k], B_n).
\]
But this limit is $0$ because the image of the map
\[
\Hom_{\T}(E_{\alpha}[k], B_{n-1}) \to \Hom_{\T}(E_{\alpha}[k], B_n)
\]
is formed by maps that factor through $\mu$ in triangle (\ref{trinum}) and
two successive maps in a triangle are $0$ whence the claim.

The object $N$ belongs to $\U$ copying the argument in
Proposition \ref{tconsDR}.
\end{proof}

\begin{rem} Observe that we cannot simplify and take $\oplus E_i$ as the
generator of the aisle because this object does not need to be
compact when the family is not finite.
\end{rem}

\begin{prop} \label{aico} Keeping the notations as in the previous Theorem,
the associated functors $\tau^{\geq n}_\U, \tau^{\leq n}_\U :\T \to \T$
preserve coproducts with $n \in \ZZ$.
\end{prop}

\begin{proof} It is obviously enough to treat the case 
$\tau^{\geq 1}_\U$. But then it is enough to see that its essential image,
$\U^{\perp}$, is closed under the formation of coproducts in $\T$. Let
$\{G_{i} / i \in I \}$ a family of objects in $\U^{\perp}$. By
using Lemma \ref{lete} once again we have to check that:
\[ \Hom_{\T}(E_{\alpha}[k], \oplus_{i \in I}G_i) = 0
\]for every $\alpha \in A$ and $k \geq 0$. But each $E_{\alpha}[k]$ is
compact and therefore,
\[ \Hom_{\T}(E_{\alpha}[k], \oplus_{i \in I}G_i) = 
\oplus_{i \in I}\Hom_{\T}(E_{\alpha}[k], G_i) = 0
\] as wanted.
\end{proof}

\begin{rem} This result is an adaptation of Neeman-Ravenel's argument for
compactly generated localizing subcategories to aisles (See
\cite[Proposition 1.9]{Ntty}).
\end{rem}

\end{document}